

\baselineskip=14pt
\parskip=10pt
\def\halmos{\hbox{\vrule height0.15cm width0.01cm\vbox{\hrule height
  0.01cm width0.2cm \vskip0.15cm \hrule height 0.01cm width0.2cm}\vrule
  height0.15cm width 0.01cm}}
\font\eightrm=cmr8 
\font\eighttt=cmtt8
\magnification=\magstephalf

\def\1{{\overline{1}}}
\def\2{{\overline{2}}}
\parindent=0pt
\overfullrule=0in
\def\Tilde{\char126\relax}
\def\frac#1#2{{#1 \over #2}}
\bf
\centerline
{
The Number of Same-Sex Marriages in a Perfectly Bisexual Population is Asymptotically Normal
}
\rm
\bigskip
\centerline{ {\it
Shalosh B. 
EKHAD}\footnote{$^1$}
{\eightrm  \raggedright
c/o D. Zeilberger, Department of Mathematics, Rutgers University (New Brunswick),
Hill Center-Busch Campus, 110 Frelinghuysen Rd., Piscataway,
NJ 08854-8019, USA.
{\eighttt c/o zeilberg  at math dot rutgers dot edu} ,
\hfill \break
{\eighttt http://www.math.rutgers.edu/\Tilde zeilberg/} .
First version: April 27, 2011.
\hfill \break
Exclusively published in the Personal Journal of Shalosh B. Ekhad and Doron Zeilberger
\hfill \break
{\eighttt http://www.math.rutgers.edu/\Tilde zeilberg/pj.html}  .
\hfill \break
Accompanied by Maple package 
{\eighttt SameSexMarriages}
downloadable from 
\hfill \break
{\eighttt http://www.math.rutgers.edu/\Tilde zeilberg/tokhniot/SameSexMarriages} .
\hfill \break
More (computer-generated!) detailed output can be gotten from the front of this article
\hfill \break
{\eighttt http://www.math.rutgers.edu/\Tilde/mamarim/mamarimhtml/ssm.html} .
}
}
{\bf Theorem:}
Consider  a population of $2n$ men and $2n$ women  where every individual is equally attracted to either sex
and chooses his or her mate according to other criteria.
Also assume that everyone gets married.
Then the expectation of the random variable ``Number of same-sex marriages" is
$$
{\frac {2n \left(2\,n  -1 \right) }{4\,n -1}}
$$
that asymptotically is
$$
n-\frac{1}{4}-\frac{1}{16}\,{n}^{-1}-{\frac {1}{64}}\,{n}^{-2}-{\frac {1}{256}}\,{n}^{-
3}-{\frac {1}{1024}}\,{n}^{-4}-{\frac {1}{4096}}\,{n}^{-5}-{\frac {1}{
16384}}\,{n}^{-6}-{\frac {1}{65536}}\,{n}^{-7}-{\frac {1}{262144}}\,{n}
^{-8}+O \left( {n}^{-9} \right)  \quad ,
$$
and its variance is
$$
8\,{\frac {{n}^{2} \left( 1-4\,n+4\,{n}^{2} \right) }{-3+28\,n-80\,{n}^
{2}+64\,{n}^{3}}} \quad,
$$
that asymptotically is
$$
\frac{1}{2}\,n+\frac{1}{8}+\frac{1}{16}\,{n}^{-1}+{\frac {3}{64}}\,{n}^{-2}+{\frac {19}{512}}\,
{n}^{-3}+{\frac {59}{2048}}\,{n}^{-4}+{\frac {45}{2048}}\,{n}^{-5}+{
\frac {17}{1024}}\,{n}^{-6}+{\frac {1637}{131072}}\,{n}^{-7}+{\frac {
4917}{524288}}\,{n}^{-8}+O \left( {n}^{-9} \right)  \quad .
$$
Furthermore, this random variable is asymptotically normal.

{\bf Semi-Rigorous Proof:}
                 The probability generating function is (why?)
$$
P_n(x) \, = \, 
\frac{1}{(4n)!/((2n)!2^{2n})}
\sum _{k=0}^{n} {{2n} \choose {2k}}^2 (2n-2k)! ((2k)!/(k!2^k))^2 x^{2k}
\quad .
\eqno(PGF)
$$
By repeatedly applying the operation $f(x) \rightarrow xf'(x)$,
plugging in $x=1$, for $n=1$ to $n=200$,
and fitting the data by a rational function of $n$, 
one gets in turn the mean, variance, and higher moments,
from which one easily gets the {\it moments about the mean} and
the {\it normalized moments} ($\alpha$ coefficients).

It turned out that the mean and variance are indeed as stated by
the theorem. As for the higher moments, we will show that the first $14$  normalized moments
tend (as $n \rightarrow \infty$) to those of the normal distribution, and the readers can go on as far as they please.

The normalized third moment (about the mean) is
$$
1/2\,\sqrt {2}\sqrt {{\frac {-3+4\,n}{{n}^{2} \left( 25-140\,n+276\,{n}
^{2}-224\,{n}^{3}+64\,{n}^{4} \right) }}} \quad ,
$$
the asymptotics of its square is:
$$
1/32\,{n}^{-5}+{\frac {11}{128}}\,{n}^{-6}+{\frac {85}{512}}\,{n}^{-7}+
{\frac {571}{2048}}\,{n}^{-8}+{\frac {3569}{8192}}\,{n}^{-9}+O \left( {
n}^{-10} \right)  \quad .
$$
Note that it goes to 0, as it should, since the odd moments of the normal distribution are all 0.

The normalized fourth moment about the mean (alias {\it kurtosis}) is
$$
1/2\,{\frac {3-100\,n+650\,{n}^{2}-1896\,{n}^{3}+2632\,{n}^{4}-1664\,{n
}^{5}+384\,{n}^{6}}{{n}^{2} \left( 35-188\,n+348\,{n}^{2}-256\,{n}^{3}+
64\,{n}^{4} \right) }} \quad .
$$
Its asymptotics is:
$$
3-{n}^{-1}+1/4\,{n}^{-2}+{\frac {7}{16}}\,{n}^{-3}+{\frac {57}{64}}\,{n
}^{-4}+{\frac {431}{256}}\,{n}^{-5}+{\frac {3137}{1024}}\,{n}^{-6}+{
\frac {22431}{4096}}\,{n}^{-7}+{\frac {159025}{16384}}\,{n}^{-8}+{
\frac {1122319}{65536}}\,{n}^{-9}+O \left( {n}^{-10} \right)  \quad .
$$
Note that it converges to $4!/(2^2 \cdot 2!)=3$, 
as it should, this being the fourth moment of the normal distribution.

The normalized fifth moment about the mean can be viewed in
\hfill\break
{\eighttt http://www.math.rutgers.edu/\Tilde zeilberg/tokhniot/oSameSexMarriages1}, and
the asymptotics of its square is:
$$
{\frac {25}{8}}\,{n}^{-5}+{\frac {315}{32}}\,{n}^{-6}+{\frac {3501}{128
}}\,{n}^{-7}+{\frac {37067}{512}}\,{n}^{-8}+{\frac {383273}{2048}}\,{n}
^{-9}+O \left( {n}^{-10} \right)  \quad .
$$
Note that it goes to 0, as it should, since the odd moments of the normal distribution are all $0$.

The normalized sixth moment about the mean can be viewed in
\hfill\break
{\eighttt http://www.math.rutgers.edu/\Tilde zeilberg/tokhniot/oSameSexMarriages1} .
Its asymptotics is:
$$
15-15\,{n}^{-1}+{\frac {31}{4}}\,{n}^{-2}+{\frac {73}{16}}\,{n}^{-3}+{
\frac {839}{64}}\,{n}^{-4}+{\frac {8401}{256}}\,{n}^{-5}+{\frac {86191}
{1024}}\,{n}^{-6}+{\frac {903617}{4096}}\,{n}^{-7}+{\frac {9635359}{
16384}}\,{n}^{-8}+{\frac {103978545}{65536}}\,{n}^{-9}+O \left( {n}^{-
10} \right)  \quad .
$$
Note that it converges to $6!/(2^3 \cdot 3!)=15$, as it should,
this being the  sixth moment of the normal distribution.

Etc. etc.  (See {\eighttt http://www.math.rutgers.edu/\Tilde zeilberg/tokhniot/oSameSexMarriages1} for evidence up to the
$14^{th}$ moment, and you are welcome to modify 
{\eighttt http://www.math.rutgers.edu/\Tilde zeilberg/tokhniot/inSameSexMarriages1} and run it on your computer,
if you want more evidence.) \halmos

{\bf Comments}

{\bf 1.} It should be routine, once we found 
the average and variance as above (that even humans should be able to do easily),
to derive a {\it local limit law}, by using Stirling's asymptotic formula.
However, for some reason, Maple refuses to take the appropriate limit.
(It works for the binomial distribution).

{\bf 2.}  Another {\it automatic}, this time 
{\it fully rigorous} proof, can also be gotten using an extension of the methods
in [Z1], using the Maple package
\hfill\break
{\tt http://www.math.rutgers.edu/\Tilde zeilberg/tokhniot/CLT}. However, since the polynomials
$P_n(x)$ are not closed-form, but rather satisfy a second-order linear recurrence (easily derived
via the Zeilberger algorithm) one should be able to do it (but {\tt CLT} has to be slightly extended).

{\bf 3.} The present approach follows [Z2] 
and the Maple package accompanying it:
\hfill\break
{\tt http://www.math.rutgers.edu/\Tilde zeilberg/tokhniot/HISTABRUT}.

{\bf 4.} Not directly related to the present article, but nevertheless  highly recommended, are
references [AZ1], [Z3], [Z4], and [Z5].

{\bf References}

[AZ1] Moa Apagodu and Doron Zeilberger,
{\it Multi-variable Zeilberger and Almkvist-Zeilberger algorithms and the sharpening of Wilf-Zeilberger theory},
Advances in Applied Mathematics {\bf 37} (2006) 139-152.

[Z1] Doron Zeilberger,
{\it The Automatic Central Limit Theorems Generator (and Much More!)},
``Advances in Combinatorial Mathematics: Proceedings of the Waterloo Workshop in Computer Algebra 2008 in honor of Georgy P. Egorychev'', 
chapter 8, pp. 165-174, (I.Kotsireas, E.Zima, eds. Springer Verlag, 2009.) Available on-line from
\hfill\break
{\tt http://www.math.rutgers.edu/~zeilberg/mamarim/mamarimhtml/georgy.html} .

[Z2] Doron Zeilberger,
{\it HISTABRUT: A Maple Package for Symbol-Crunching in Probability theory},
The Personal Journal of Shalosh B. Ekhad and Doron Zeilberger, Aug. 25, 2010.
(Journal Main page: {\tt http://www.math.rutgers.edu/~zeilberg/pj.html}, direct url of article:
\hfill\break
{\tt http://www.math.rutgers.edu/~zeilberg/mamarim/mamarimhtml/histabrut.html})

[Z3] Doron Zeilberger,
{\it Andr\'e's  reflection proof generalized to the many-candidate ballot problem},
Discrete Math {\bf 44} (1983), 325-326.

[Z4] Doron Zeilberger,
{\it Kathy O'Hara's constructive proof of the unimodality of the Gaussian polynomials},
Amer. Math. Monthly {\bf 96}(1989), 590-602.

[Z5] Doron Zeilberger,
{\it The Umbral Transfer-Matrix Method. I. Foundations}
Journal of Combinatorial Theory, Series A,  {\bf 91} (2000), 451-463.

\end